\documentclass[12pt,reqno]{amsart}

\usepackage{latexsym,color,amsmath,amsthm,amssymb,amscd,amsfonts}
\usepackage{tikz}
\usetikzlibrary{arrows}
\usepackage{scalefnt}
\usepackage[font=small,labelfont=bf]{caption}
\usepackage{float}
\usepackage{graphicx}
\usepackage{relsize}
\usepackage{amsopn}
\usepackage{mathrsfs}  
\usepackage[hidelinks]{hyperref}
\usepackage{bm}
\usepackage{bbm}
\usepackage[shortlabels]{enumitem}
\usepackage{multicol}
\usepackage{float} 
\usepackage{todonotes}
\usepackage{cite}

\usepackage[margin=2.5cm]{geometry}

\newtheoremstyle{nonum}{}{}{\itshape}{}{\bfseries}{.}{ }{\thmnote{#3}}

\newtheorem{thm}{Theorem}[section]
\newtheorem*{thm*}{Theorem}

\newtheorem{cor}[thm]{Corollary}
\newtheorem{lem}[thm]{Lemma}

\newtheorem{conj}{Conjecture}

\newtheorem*{definition*}{Definition}
\newtheorem{fact}[thm]{Fact}

\newtheorem*{rems*}{Remarks}
\newtheorem{rems}[thm]{Remarks}

\theoremstyle{nonum}

\newcommand{\R}{\mathbb R}
\newcommand{\RR}{\mathbb R}
\newcommand{\RRn}{\RR_+}

\newcommand{\Z}{\mathbb Z}

\def\sp{ {\rm{sp} }  }

\def\vol{{\rm Vol}}

\def\sp{\spade}

\def\conv{{\rm conv}}

\newtheoremstyle{nonum}{}{}{\itshape}{}{\bfseries}{.}{ }{\thmnote{#3}}

\newtheorem*{lem*}{Lemma}

\theoremstyle{remark}

\def\sp{{\rm sp }}

\def\conv{{\rm conv}}

\keywords{locally anti-blocking bodies, simplices, (mixed) volume
inequalities, Godbersen's conjecture, mixed volume formulas}


\begin{document}

\title {Godbersen's conjecture for locally anti-blocking bodies}

\author{Shay Sadovsky}
\address{School of Mathematical Sciences, Tel Aviv University, Israel}
\email{shaysadovsky@mail.tau.ac.il}

\maketitle
\begin{abstract}
	In this note we give a short proof of Godbersen's conjecture for the class of locally anti-blocking bodies. We show that all equality cases amongst locally anti-blocking bodies are for simplices, further supporting the conjecture. The proof of equality cases introduces a useful calculation of mixed volumes of aligned simplices.
\end{abstract}

\section{Introduction and results}\label{sec:intro}
\newcommand\Def[1]{\textbf{#1}}%
\newcommand\PM{\{-1,1\}}%

In this note we will prove Godbersen's conjecture \cite{Godbersen} for a special class of sets called locally anti-blocking
bodies. 
\begin{conj}
	For any convex body $K \subseteq\RR^n$ and $0 < j < n$,
	\begin{equation}\label{eqn:Godberson}
	V_n(K[j],-K[n-j]) \ \leq \ \binom{n}{j} \vol (K).	 
	\end{equation}
	Equality holds for $j\neq 0, n$ if and only if $K$ is a simplex.
\end{conj}
This conjecture (also conjectured independently in \cite{Makai}) is a refinement to Rogers and Shephard's inequality for the difference body \cite{RS-difference}, and is known for a handful of special cases. For more background and previous results, see \cite{Artstein-Note,ArtsteinKeshet,AB-bodies}.

A convex body $K \subseteq\RR^n$ is called \Def{$1$-unconditional} if
$(x_1,\dots,x_n) \in K$ implies that the points $(\pm x_1,\dots, \pm x_n)$
also belong to $K$. Clearly such a body is determined by its subset  $K_+ := K
\cap \RRn^n$, which is itself a convex set of a special kind, called an
``anti-blocking body''~\cite{Fulkerson1971blocking, Fulkerson1972anti} or a
``convex corner''~\cite{Bollobas-Reverse-Kleitman}. A convex body $K \subset
\RRn^n$ is called \Def{locally anti-blocking} if, for any coordinate hyperplane $E^c_J = \sp \{e_j\}_{j\in J}$, $J\subset [n]$ (where $\{e_i\}_{i=1}^n$ is the standard basis), one has $P_E K = K\cap E$. Alternatively, they can be defined as bodies that are anti-blocking in each orthant. For an extended review of these classes of bodies see \cite{AB-bodies}.

Recently, the class of locally anti-blocking has recieved more attention as a natural extension of the notion of unconditional convex bodies, and several conjectures have been verified for this class -- Mahler's conjecture \cite{AB-bodies, FradeliziMeyer-func-santalo} and Kalai's $3^d$ conjecture \cite{Sanyal-3d}.

We prove the following theorem.
\begin{thm}\label{thm:god-for-LAB}
	Locally anti-blocking bodies satisfy Godbersen's conjecture.
\end{thm}

The methods in this note are based on several inequalities obtained for anti-blocking bodies in \cite{AB-bodies}, along with new geometric observations on anti-blocking bodies and simplices.

The proof of the equality cases requires a useful geometric lemma regarding mixed volumes of ``coordinate-aligned" simplices, originally due to Shephard \cite{shephard1960} and extended in \cite{henk2012steiner} by Henk, Hern\'andez Cifre and Saor\'in. We denote by $\Delta_n$ the \emph{coordinate simplex} $\conv\{0^n, e_1,\dots, e_n\}\in \RR^n$. For completeness, we include two proofs of this lemma in section \ref{sec:simplex}.
\begin{lem}\label{lem:mixed-of-simpex}
	Let $\alpha_1,\dots,\alpha_n\ge0$ and let  $K=\conv(0,\alpha_1 e_1,\dots, \alpha_n e_n)$. Then, for any $0 \le j \le n$,
	\[V_n(K[j],\Delta_n[n-j])=\frac{1}{n!} \max\left \{\prod_{i\in I} \alpha_i: |I|=j \right \}.\]
\end{lem}

\subsection*{Acknowledgments}
The author would like to thank Shiri Artstein-Avidan, for her supervision and guidance, and Arnon Chor, for many insightful discussions and his thorough reading of this manuscript. The author thanks Raman Sanyal for pointing out the second proof of Lemma \ref{lem:mixed-of-simpex}, and Martin Henk for the bringing to our attention its existence in \cite{shephard1960, henk2012steiner}. This research was partially supported by ISF Grant No. 784/20. The author is also grateful to the Azrieli foundation for the award of an Azrieli fellowship.

\section{Definitions and previous results for locally anti-blocking bodies}\label{sec:Def}
In this section we present several lemmas describing properties of locally anti-blocking bodies and which will be used in the proof of the theorem. We begin with some definitions.

For $\sigma \in \PM^n$ and $S \subseteq\R^n$ let us write $\sigma S = \{ (
\sigma_1 x_1, \dots, \sigma_n x_n) : x \in S \}$. A convex body $K
\subseteq\RR^n$ is locally anti-blocking if $(\sigma K) \cap \RRn^n$
is anti-blocking for all $\sigma \in \PM^n$. In particular, the distinct bodies
$K_\sigma := K \cap \sigma \RRn^n$ form a dissection of $K$.  That is,
$K = \bigcup_\sigma K_\sigma$ and the $K_\sigma$ have disjoint interiors. By
definition, $1$-unconditional convex bodies are precisely the locally
anti-blocking bodies with $\sigma K_\sigma = K_+$ for all $\sigma$.

The \Def{mixed volume}, introduced by Minkowski,
is the non-negative symmetric function $V_n$ defined on $n$-tuples of convex
bodies in $\RR^n$  satisfying
\begin{equation}\label{eq:mixvol-formula}
\vol_n(\lambda_1 K_1 + \dots + \lambda_n K_n) = \sum_{i_1,\dots, i_n =1}^n
\lambda_{i_1} \cdots \lambda_{i_n} V_n(K_{i_1},\dots, K_{i_n}),
\end{equation}
for any convex bodies $K_1,\dots, K_n \subset \RR^n$ and $\lambda_1, \dots,
\lambda_n \ge 0$; see~\cite{Schneider-book-NEW} for details.  We call
$V_n(K_1,\dots, K_n)$ the mixed volume of $K_1,\dots, K_n$. We  denote the mixed volume of $j$ copies of $K$ and $n-j$ copies of  $T$ by $V_n(K[j], T[n-j])$.

The following lemma from \cite{AB-bodies} (see also \cite{Sanyal-ABbodies}) gives a formula for the mixed volume of anti-blocking bodies in opposite orthants.
\begin{lem}\label{lem:AB-decomposition-mixed}
	Let $K, K' \subseteq \RRn^n$ be anti-blocking bodies, let $0\le j\le n$. Then
	\begin{align}\label{eq:ABjmixformula}
	V_n(K[j], -K'[n-j]) \ &= \ {\binom{n}{j}}^{-1}\sum_{\substack{E \text{ a }j \text{-dim.}\\ \text{coord. hyperplane}}}
	\vol_j(P_E K) \cdot \vol_{n-j} (P_{E^\perp}K') \\
	\vol(K \vee -K') \ &= \  \sum_{j=0}^n   V_n(K[n-j],-K'[j])
	\label{eq:vol-j-conv(K,-t)} \, .
	\end{align}
\end{lem}

 It turns out that the mixed volume of two locally anti-blocking bodies decomposes into the mixed volumes in each orthant. 
\begin{lem}
	Let $K, K' \subseteq \RRn^n$ be locally anti-blocking bodies. Then,
	\begin{equation}\label{eq:LAB-vol-decomp}
	\vol(K+K')=\sum_{\sigma\in \{-1,1\}^n } \vol(K_\sigma+K'_\sigma)
	\end{equation}
	and in fact,
	\begin{equation}\label{eq:LAB-j-mixvol}
		V_n(K[j], K'[n-j]) = \sum_{\sigma\in \{-1,1\}^n } V_n(K_\sigma[j], K'_\sigma[n-j]).
	\end{equation}
\end{lem}

\begin{proof}
	In \cite[Lemma 2.2]{AB-bodies} it was shown that 
	$$ K + K'  \ = \ \bigcup_\sigma K_\sigma + K'_\sigma \ . $$
	Since the bodies $K_\sigma + K'_\sigma$ are in different orthants, they are disjoint up to measure zero, implying \eqref{eq:LAB-vol-decomp}.
	
	Note that for $\lambda \ge 0$, the body $\lambda K$ is also locally anti-blocking. Expanding both sides of \eqref{eq:LAB-vol-decomp} using the formula for mixed volume given by \eqref{eq:mixvol-formula}, we get that
	\[\sum_{j=1}^n \lambda^j \binom{n}{j} V_n(K[j],K'[n-j]) = \sum_{\sigma\in \{-1,1\}^n} \lambda^j \binom{n}{j} V_n(K_\sigma[j],K'_\sigma[n-j]).\]
	Comparing the coefficients of the polynomial in $\lambda$ on both sides proves \eqref{eq:LAB-j-mixvol}.
\end{proof}

We will also make use of the following simple observation for locally anti-blocking bodies.
\begin{fact}\label{fact:subspaces-proj-equality}
	Let $K\subset \RR^n$ be locally anti-blocking, and let $E:=\sp\{e_i: i\in I\}\subset \RR^n$ for some $I\subset [n]$. Let $\tau, \sigma \in \PM^n$ be two orthant signs, such that $\tau|_I = \sigma |_I$. Then,
	\[P_E K_\sigma = P_E K_\tau.\]
\end{fact}
\begin{proof}
	Since $\tau$ and $\sigma$ are equal in the coordinates of $E$, $$\tau \RR_+^n \cap E = \sigma \RR_+^n \cap E.$$
	Recalling the definition of anti-blocking bodies, 
	\[P_E K_\sigma = K_\sigma \cap E = K \cap \sigma \RR_+^n \cap E = K \cap \tau \RR_+^n \cap E = P_E K_\tau .\]
\end{proof}

The last lemma we recall is a mixed volume version of the Reverse Kleitman inequality of
Bollob\'as, Leader, and Radcliffe~\cite{Bollobas-Reverse-Kleitman} which was proved in \cite{AB-bodies} as well.

\begin{thm}\label{thm:kleitman-bound-on-mixed}
	Given two anti-blocking bodies, $K, T \subseteq\RR_+^n$,
	\begin{equation}\label{eq:reverse-kleitman}
	V_n(K[j],T[n-j]) \ \leq \  V_n(K[j],-T[n-j]).
	\end{equation}
	In particular, 
	\[ \vol(K+T) \ \le \ \vol(K-T). \] 
\end{thm}

\section{Proof of the main theorem}\label{sec:Pf}
In this section we prove Theorem \ref{thm:god-for-LAB}. The proof is in two parts, first we prove the inequality, and then show that all equality cases are simplices.
\begin{proof}[Proof of the inequality in Theorem \ref{thm:god-for-LAB}]
Let $K=\cup_\sigma K_\sigma\subset \RR^n$ be locally anti-blocking. Note that for any $\sigma \in \PM$, $(-K)_\sigma =-( K_{-\sigma})$.

Applying \eqref{eq:LAB-j-mixvol} and \eqref{eq:reverse-kleitman} in the orthant $\sigma \RR_+^n$, we get
\begin{align}
V_n(K[j],-K[n-j])&=\sum_{\sigma \in \PM^n} V_n(K_\sigma[j],(-K)_\sigma [n-j]) \nonumber \\
&=\sum_{\sigma \in \PM^n} V_n(K_\sigma[j],-(K_{-\sigma}) [n-j]) \nonumber \\
&\le \sum_{\sigma \in \PM^n} V_n(K_\sigma[j],(K_{-\sigma}) [n-j]). \label{eq:first-ineq}
\end{align}
Noting that $K_\sigma$ and $-(K_{-\sigma})$ are anti-blocking and in the same orthant, we may apply \eqref{eq:ABjmixformula} to each mixed volume in the above sum and get
\begin{align*}
 \sum_{\sigma \in \PM^n} V_n(K_\sigma[j],(K_{-\sigma}) [n-j]) = \sum_{\sigma \in \PM^n} \binom{n}{j}^{-1}  \sum_{\substack{E \text{ a }j \text{-dim.}\\ \text{coord. hyperplane}}} \vol_j(P_E K_\sigma)\vol_{n-j}(P_{E^\perp} K_{-\sigma}).
\end{align*}
Fixing $E=\sp\{e_i\}_{i\in I}$ a $j$-dimensional coordinate subspace, let $\sigma\in \PM^n$, we claim that there is a unique $\tau \in \PM^n$ such that 
\[P_EK_\sigma = P_EK_\tau \text{ and } P_{E^\perp}K_{-\sigma} = P_{E^\perp} K_\tau\]
and that this is a bijection. The map $\tau \mapsto \sigma(\tau, E)$, $\tau \in \PM^n$, where $\sigma=\sigma(\tau,E) \in \PM^n$ is the unique vector satisfying
\[\sigma|_I =\tau|_I,\ \ \sigma|_{I^c}=-\tau|_{I^c}\] is a bijection (indeed, it is clearly inversed by itself).
By Fact \ref{fact:subspaces-proj-equality} we have $P_E K_\tau = P_E K_{\sigma(\tau, E)}$ by the definition of $\sigma(\tau, E)$, and similarly $P_{E^\perp} K_\tau = P_{E^\perp} K_{-\sigma(\tau, E)}$.
We may now apply this bijection and get that 
\begin{align*}
 &\sum_{\sigma \in \PM^n} \vol_j(P_E K_\sigma)\vol_{n-j}(P_{E^\perp} K_{-\sigma}) \\
 &=  \sum_{\tau \in \PM^n} \vol_j(P_E K_{\sigma(\tau, E)})\vol_{n-j}(P_{E^\perp} K_{-\sigma(\tau, E)})\\
 & = \sum_{\tau \in \PM^n} \vol_j(P_{E} K_\tau )\vol_{n-j}(P_{E^\perp} K_{\tau}),
\end{align*}

Changing the order of summation and applying the classical Rogers-Shephard inequality for sections and projections \cite{RS-sections-projections} (recalling that for locally anti-blocking bodies $K_\tau \cap E = P_E K_\tau$) we get
\begin{align} 
&\sum_{\sigma \in \PM^n} \binom{n}{j}^{-1}  \sum_{\substack{E \text{ a }j \text{-dim.}\\ \text{coord. hyperplane}}} \vol_j(P_E K_\sigma)\vol_{n-j}(P_{E^\perp} K_{-\sigma}) \nonumber\\
&=\binom{n}{j}^{-1}  \sum_{\substack{E \text{ a }j \text{-dim.}\\ \text{coord. hyperplane}}} \sum_{\sigma \in \PM^n} \vol_j(P_E K_\sigma)\vol_{n-j}(P_{E^\perp} K_{-\sigma}) \nonumber\\
&=\binom{n}{j}^{-1}  \sum_{\substack{E \text{ a }j \text{-dim.}\\ \text{coord. hyperplane}}}\sum_{\tau \in \PM^n} \vol_j(P_E K_\tau )\vol_{n-j}(P_{E^\perp} K_{\tau})\nonumber \\
 &\le \binom{n}{j}^{-1} \sum_{\substack{E \text{ a }j \text{-dim.}\\ \text{coord. hyperplane}}}  \sum_{\tau \in \PM^n} \binom{n}{j}\vol_n(K_\tau) \label{eq:second-ineq} \\
 & = \sum_{\substack{E \text{ a }j \text{-dim.}\\ \text{coord. hyperplane}}} \sum_{\tau \in \PM^n}\vol_n(K_\tau) = \binom{n}{j} \vol(K) \nonumber
\end{align}
completing the proof.

\end{proof}

\begin{proof}[Proof of the equality cases in Theorem \ref{thm:god-for-LAB}]
We show that the only locally anti-blocking body $K$ for which $V_n(K[j],-K[n-j])=\binom{n}{j}\vol(K)$ is a simplex. The only inequalities in the proof are in \eqref{eq:first-ineq} and \eqref{eq:second-ineq}.

Equality in \eqref{eq:second-ineq} implies that for every $j$-dimensional coordinate subspace, 
\[\vol_j(P_E K_\sigma )\vol_{n-j}(P_E K_{\sigma}) = \binom{n}{j}\vol_n(K_\sigma).\]
It was shown in \cite[Proposition 3.3]{AB-bodies} that for anti-blocking bodies of full dimension, this implies that $K_\sigma$ is a simplex. Thus, we can conclude that $K$ is of the form $\cup_{\sigma \in \{-1,1\}^n} K_\sigma$ where $K_\sigma = \conv\{\alpha_{i,\sigma} \sigma_i e_i\}_{i=1}^n$ and $\alpha_{i,\sigma}\ge 0$.
It turns out, and will follow from the proof below, that the only simplices which can be constructed in this manner are, up to change of coordinates, of the form
\begin{equation} \label{eq:equality-cases}
K=\conv\left(\{0_n\} \cup \{\alpha_i  e_i\}_{i=1}^n\right) \text{ or }K=\conv\left( \{\alpha_1 e_1, -\beta_1 e_1\}\cup \{\alpha_i  e_i\}_{i=2}^n \right)
\end{equation}
for some $\alpha_i>0, \beta_1>0$.

Assume, without loss of generality, that $K_{(1,\dots,1)}=\Delta_n$, and $K_{-1} =K_{(-1,\dots,-1)}=\conv(\{-\alpha_i e_i\}_{i=1}^k\cup \{0_n\})$ for some $\alpha_i\ge \alpha_{i+1}>0$ and $0\le k\le n$ (indeed, one orthant must be of full dimension, and we may apply a linear transformation on $K$ so that it is the first orthant and so that the simplex in this orthant is $\Delta_n$). 
If $k=0$ then there is a trivial equality in \eqref{eq:first-ineq}, which corresponds to the case of $$K=\conv\left(\{0_n\} \cup \{\alpha_i  e_i\}_{i=1}^n\right).$$

According to Fact \ref{fact:simplex-sum-lower-dim} from the next section, 
\begin{align*}
	&\vol(\Delta_n - \lambda K_{-1})=\int_{(t_1,\dots, t_{n-k})\in \Delta_{n-k}} \vol_k((1-\sum_{i=1}^{n-k} t_i)\Delta_k - \lambda K_{-1}) dt\\
	&=\sum_{j=0}^k\lambda^{k-j}\binom{k}{j}V_k(\Delta_k [j],-K_{-1}[k-j])  \int_{(t_1,\dots, t_{n-k})}(1-\sum_{i=1}^{n-k} t_i)^j dt .
\end{align*}
and similarly
\[\vol(\Delta+ \lambda K_{-1}) = \sum_{j=0}^k\lambda^{k-j}\binom{k}{j}V_k(\Delta_k [j],K_{-1}[k-j])  \int_{(t_1,\dots, t_{n-k})}(1-\sum_{i=1}^{n-k} t_i)^j dt.\]
Since the integral does not depend on $K$, comparing coefficients shows that the $j$-th mixed volume in $\RR^n$ $V_n(\Delta_n[j],-K[n-j])$ is proportional to $V_k(\Delta_k[j],-K[k-j])$, and similarly for $K$, so that equality in \eqref{eq:first-ineq} holds for simplices $\Delta_n$ and $K$ if and only if equality holds for the mixed volumes in $\RR^k$, i.e. if
\begin{equation}\label{eq:equality-cond}
V_k(\Delta_k [j],-K_{-1}[k-j])= V_k(\Delta_k [j],K_{-1}[k-j]).
\end{equation}
If $k=1$ then the equality above holds trivially, and $K_{-1} = \conv(-\alpha_1 e_1,0)$. This clearly implies that
$$K=\conv\left( \{ e_1, -\alpha_1 e_1\}\cup \{e_i\}_{i=2}^n \right),$$ which is, up to a linear transformation, the second case in \eqref{eq:equality-cases} and is also a simplex.

Let us show that $k=0,1$ are the only equality cases, i.e. that  \eqref{eq:equality-cond} does not hold for $k\ge 2$. Note that $\Delta_k$ and $K_{-1}$ are both simplices, and they satisfy the assumptions of Lemmas \ref{lem:mixed-of-simpex} and \ref{lem:AB-decomposition-mixed}. Using these Lemmas we compute
\begin{align*}
V_k(\Delta_k [j],-K_{-1}[k-j])&=\frac{1}{k!}\prod_{i=1}^j \alpha_i\\
V_k(\Delta_k [j],K_{-1}[k-j])&=\binom{k}{j}^{-1}\sum_{\substack{E \text{ a }j \text{-dim.}\\ \text{coord. hyperplane}}} \vol_j(P_E \Delta_k)\vol_{k-j}(P_{E^\perp} K_{-1})\\
&=\binom{k}{j}^{-1} \sum_{J\subset [k],|J|=j} \frac{\prod_{i\in J}\alpha_i}{j!(k-j)!}= \frac{1}{k!}\sum_{J\subset [k],|J|=j} {\prod_{i\in J}\alpha_i}
\end{align*}
and since $\alpha_i>0$ for all $i=1,\dots k$ and $j\ge 1$, \eqref{eq:equality-cond} does not hold, and so there is an inequality in \eqref{eq:first-ineq} and subsequently in Theorem \ref{thm:god-for-LAB}.

We have thus shown that amongst all locally-anti-blocking bodies, the equality cases in Theorem \ref{thm:god-for-LAB} are only simplices.
\end{proof}
\section{Mixed volumes of simplices}\label{sec:simplex}
In this section, we prove Lemma \ref{lem:mixed-of-simpex} and make some additional observations regarding mixed volumes of two simplices. We start with the following simple fact, which was already used in the previous section and which is also used to prove the lemma.

\begin{fact}\label{fact:simplex-sum-lower-dim}
	Let $E\subseteq \RR^n$ be the $k$-dimensional coordinate hyperplane given by $E=\sp\{e_n,\dots, e_{n-k+1}\}$ and let $K\subseteq E$ then, 
	\[\vol_n(\Delta_n + K) =\int_{(t_1, \dots,t_{n-k})\in \Delta_{n-k}} \vol_{k}\left(  (1-\sum_{i-1}^{n-k}t_i)\Delta_k + K \right) dt.\]
\end{fact}
\begin{proof}
	Due to Fubini's theorem,
	\[\vol_n(\Delta_n + K) = \int_{t \in E^\perp \cap \Delta_n + K} \vol_k( (\Delta_n +K)\cap (t+E) ) dt.\]
	Since $\Delta_n= \{x:\ \sum_{i=1}^n x_i \le 1, \ 0\le x_i\}$, we get that 
	\begin{align*}
	&(\Delta_n +K)\cap (t+E) =  \\
	&= \{(t_1,\dots, t_{n-k}, x_{n-k +1}+y_{n-k+1},\dots,x_n+y_n): \ \sum_{i=n-k_+1}^n x_i \le (1-\sum_{i=1}^{n-k}t_i),\  x_i\ge 0,\  y\in K\}\\
	&= t+ ((1-\sum_{i=1}^{n-k}t_i))\Delta_E +K.
	\end{align*}
	So, the volume of this set is $\vol_{k}\left(  (1-\sum_{i=1}^{n-k}t_i)\Delta_k + K \right)$, completing the proof.
\end{proof}

We present two different proofs of Lemma \ref{lem:mixed-of-simpex}. The first is an elementary, geometric proof, in the spirit of the original proof in \cite{shephard1960}. The second is an algebraic proof via the celebrated Bernstein--Khovanskii--Kouchnirenko (BKK) Theorem \cite{Bernstein, Khovanskii, Kouchnirenko}, which states the following:

\begin{thm}[Bernstein--Khovanskii--Kouchnirenko (BKK) \cite{Bernstein, Khovanskii, Kouchnirenko}]
	Given polynomials
	$f_1, \dots, f_n\in \mathbb{C}[x_1^{\pm 1},\dots, x_n^{\pm 1}]$, let $P_i = NP(f_i)$ be the Newton polytope of $f_i$ in $\RR^n$. Then, for generic choices of the coefficients in the
	$f_i$, the number of common solutions (with multiplicity) is exactly $V_n(P_1, . . . , P_n)$.
\end{thm}

Both proofs are presented below.

\begin{proof}[First (elementary) proof of Lemma \ref{lem:mixed-of-simpex}]
	Assume without loss of generality that $\alpha_i \ge \alpha_{i+1}$ for $1\le i \le n$, since applying a linear transformation that changes the order of the coordinates does not change $\Delta_n$. We prove the lemma by induction on the dimension.
	
	Let $n=1$, $K\subset \RR$ as above, i.e. $\Delta_n=[0,1]$ and $K=[0,\alpha_1]$, then $j\in \{0,1\}$ and the equality holds trivially.
	
	Let $n\ge 2$, and assume we have proved the Lemma in $\RR^{n-1}$ for all $\alpha_1 \ge \dots\ge  \alpha_{n-1}\ge 0$. We first notice that 
	\begin{equation}\label{eq:sum-simplex-decompose}
	\Delta_n +K  = (e_n + K) \cup (\Delta_n+P_{e_n^\perp} K)
	\end{equation}
	and that the two parts have a zero-volume intersection. If $\alpha_n = 0$ then the equality is clear. The inclusion $\Delta + K \supseteq (e_n + K) \cup (\Delta_n+P_{e_n^\perp} K)$ is trivial and so is the zero-volume intersection, so need only to show the reverse inclusion.
	Indeed, by the definition of $K$ and $\Delta_n$, and as convex hull commutes with Minkowski addition,
	\begin{align*}
	\Delta_n+K &= \conv(\{0^n\}\cup \{e_i\}_{i=1}^n)+\conv(\{0^n\}\cup \{\alpha_i e_i\}_{i=1}^n)\\
	&=\conv(0^n \cup \{e_i + \alpha_k e_k\}_{i,k=1}^n ).
	\end{align*}
	In fact, the points of the form $e_i+\alpha_n e_n$ when $i\neq n$ can be omitted from this hull.
	Indeed, take $\lambda = \frac{\alpha_n}{1+\alpha_n}$, $\mu = \frac{1}{1+\alpha_i}$, then as $\alpha_n\le \alpha_i$, $1-\lambda-\mu \ge 1-\frac{\alpha_n}{1+\alpha_n}-\frac{1}{1+\alpha_n}=0$. Since $0^n, (1+\alpha_n)e_n, (1+\alpha_i)e_i$ are in the hull,
	\[e_i+\alpha_n e_n = \lambda (1+\alpha_n)e_n + \mu (1+\alpha_i)e_i + (1-\lambda-\mu)0^n. \]
	So we find that 
	\[ \Delta_n + K = \conv(\{0^n\} \cup \{e_i+\alpha_k e_k\}_{i=1, k=1}^{n, n-1} \cup \{e_n+\alpha_n e_n\}).\]

	Note that
	\[\Delta_n +P_{e_n^\perp} K = \conv(\{0^n\} \cup \{e_n\} \cup \{e_i+\alpha_k e_k\}_{i=1, k=1}^{n, n-1}),\]
	so that $\Delta_n +K$ can be described as
	\begin{align*}
	\Delta_n+ K &= \conv(\conv(\{0^n\} \cup \{e_n\} \cup \{e_i+\alpha_k e_k\}_{i=1, k=1}^{n, n-1}) \cup \{(1+\alpha_n) e_n\})\\
	&= \conv((\Delta_n +P_{e_n^\perp} K) \cup \{(1+\alpha_n)e_n\}).
	\end{align*}
	
	This description shows that $\Delta_n+ K$ is equal to the union of all segments connecting between $(1+\alpha_n)e_n$ and points in $\Delta_n +P_{e_n^\perp} K$. Denote by $\ell=\ell(x,y)$ such a segment from $z=(1+\alpha_n) e_n $ to a point $x\in \Delta_n +P_{e_n^\perp} K$. Let us show that $\ell$ intersects the set $e_n +P_{e_n^\perp} K $ at some point $y$ which is in $e_n +P_{e_n^\perp} K$.
	Due to the convexity of $\Delta_n +P_{e_n^\perp} K $ and $e_n+P_{e_n^\perp} K$, it will suffice to find such a point $y$ only for $x\in \{0^n\} \cup \{e_n\} \cup \{e_i+\alpha_k e_k\}_{i=1, k=1}^{n, n-1}$, as this set contains the extremal points of  $\Delta_n +P_{e_n^\perp} K$.
	
	Indeed, the line connecting $z$ with $x=(1+\alpha_k)e_k$ intersects the edge $\conv\{e_n,e_n+\alpha_k e_k\} \in e_n +P_{e_n^\perp} K$ at the point $y= e_n+\frac{\alpha_n}{1+\alpha_n}(1+\alpha_k)e_k$. In a similar fashion, the line connecting  $z$ with $x=e_i+\alpha_ke_k$ intersects the face $\conv\{e_n, e_n+\alpha_i e_i, e_n+\alpha_k e_k\}$ of $e_n +P_{e_n^\perp} K$ at $y=e_n+(\frac{\alpha_n}{\alpha_i(1+\alpha_n)})\alpha_i e_i +(\frac{\alpha_n}{(1+\alpha_n)})\alpha_k e_k$. Finally, the line between $z$ and $x=e_n$ is contained in $e_n+K$ and thus intersects $e_n +P_{e_n^\perp} K$ trivially at $e_n$, and the line between $z$ and $x=0^n$ intersects $e_n +P_{e_n^\perp} K$ at $e_n$.
	
	Thus, we have shown that $\ell(x,z)$ may be split into two segments, $\ell(z,y)\in e_n +K$ and $\ell(y,x)\in \Delta_n +P_{e_n^\perp} K$, where $y = \ell(x,z)\cap (e_n +P_{e_n^\perp} K)$. So we may prove the reverse inclusion
	\begin{align*}
	\Delta_n+ K &= \cup\{\ell(x,z): x\in \Delta_n +P_{e_n^\perp} K\} \\
	&= \cup\{\ell(x,y)\cup \ell(y,z): x\in \Delta_n +P_{e_n^\perp} K,\ y \in \ell(x,z)\cap (e_n +P_{e_n^\perp} K) \}\\
	&\subseteq \cup\{\ell(x,y): x\in \Delta_n +P_{e_n^\perp} K,\ y\in e_n +P_{e_n^\perp} K \}  \bigcup \cup\{\ell(y,z):  y\in (e_n +P_{e_n^\perp} K) \}\\
	& \subseteq (\Delta_n + P_{e_n^\perp}K) \cup (e_n +K)
	\end{align*} thus proving \eqref{eq:sum-simplex-decompose}.

	Using \eqref{eq:sum-simplex-decompose} and Fact \ref{fact:simplex-sum-lower-dim},
	\begin{align*}
	\vol( \Delta_n+\lambda K)&=\vol(e_n + \lambda K)+\vol( \Delta_n+P_{e_n^\perp}\lambda K)\\
	&=\lambda^n \vol(K)+  \int_0^1 \vol_{n-1}(t( \Delta_{n-1}) + P_{e_n^\perp} \lambda K ) dt.
	\end{align*}
	By the induction hypothesis for $\Delta_{n-1}$ and $P_{e_n^\perp } K=\conv(0,\alpha_1 e_1,\dots, \alpha_{n-1} e_{n-1})$:
	\begin{align*}
	\vol(\Delta_n+\lambda K)&=  \lambda^n \vol(K)+\int_0^1 \sum_{j=0}^{n-1} t^{n-j-1} \lambda^j \binom{n-1}{j} V_{n-1}(P_{e_n^\perp }K[j],\Delta_{n-1} [n-j-1]) dt\\
	&=\lambda^n\vol(K)+\int_0^1 \sum_{j=0}^{n-1} t^{n-j-1}\lambda^j \frac{\prod_{i=1}^{j}\alpha_i}{j!(n-j-1)!} dt\\
	& = \sum_{j=0}^{n}\lambda^j \frac{\prod_{i=1}^{j}\alpha_i}{j!(n-j)!} = \sum_{j=0}^{n}\binom{n}{j}\lambda^j \frac{\prod_{i=1}^{j}\alpha_i}{n!}.
	\end{align*}
	Comparing the coefficients of $\lambda^j$  we get
	\[V_n(\Delta_n[j],K[n-j])= \frac{\prod_{i=1}^j \alpha_j}{n!}.\]
\end{proof}

\begin{proof}[Second (algebraic) proof of Lemma \ref{lem:mixed-of-simpex}]
	Assume without loss of generality that for all $1\le i \le n$ $\alpha_i \in \Z$ and $\alpha_i \ge \alpha_{i+1}$. Hence, $\Delta_n$ is the Newton polytope of a (generic) polynomial
	\[f(x) = \sum_{i=1}^n c_i x_i + c_0, c_i \in \mathbb{C},\]
	and similarly $K$ is the Newton polytope of a (generic) polynomial
	\[g(x) = \sum_{i=1}^n c_i x_i^{\alpha_i}+c_0, c_i \in \mathbb{C}.\]
	
	According to the BKK Theorem, $V_n(K_n[j], \Delta_n [n-j])$ is equal to the number of solutions to the system of equations $g_1=\dots g_j= f_1 = \dots f_{n-j} = 0$ where $(g_i)_{i=1}^j$ are copies of $g$ with generic coefficients, and similarly $(f_i)_{i=1}^{n-j}$ are generic copies of $f$.
	
	Noting that $f_1,\dots f_{n-j}$ are linearly dependent, we may apply row operations to the system of equations until it is of the following form:
	\begin{equation*}
	\begin{cases}
	0 = c_0^\ell + x_\ell + \sum_{i=1}^j c_i^\ell x_i & j+1 \le \ell \le n \\
	g_\ell = 0 = c_0^\ell + \sum_{i=1}^j c_i^\ell x_i^{\alpha_i} + \sum_{i=j+1}^n c_i^\ell x_i^{\alpha_i} & 1\le \ell \le j. 
	\end{cases}
	\end{equation*}
	We may plug $x_\ell = -c_0^\ell - \sum_{i=1}^j c_i^\ell x_i$ into the equations of the form $g_\ell = 0$ to get
	\begin{equation*}
	\begin{cases}
	0 = c_0^\ell + x_\ell + \sum_{i=1}^j c_i^\ell x_i & j+1 \le \ell \le n \\
	g_\ell = 0 = c_0^\ell + \sum_{i=1}^j c_i^\ell x_i^{\alpha_i} + \sum_{i=j+1}^n c_i^\ell(-c_0^i -  \sum_{m=1}^j c_m^i x_m)^{\alpha_i} & 1\le \ell \le j. 
	\end{cases}
	\end{equation*}
		
	Notice that we may undo the row operations, to get back to the equivalent system
	\begin{equation*}
	\begin{cases}
	f_\ell = 0  & j+1 \le \ell \le n \\
	g_\ell = 0 = c_0^\ell + \sum_{i=1}^j c_i^\ell x_i^{\alpha_i} + \sum_{i=j+1}^n c_i^\ell(-c_0^i -  \sum_{m=1}^j c_m^i x_m)^{\alpha_i} & 1\le \ell \le j. 
	\end{cases}
	\end{equation*}	
	
	We are in position to prove the theorem. The Newton polytope of the polynomials $g_\ell$ above is of the form
	\begin{align*}
		&\conv\left( \{\alpha_i e_i\}_{i=1}^j \cup \left\{ \sum_{k=1}^j \beta_k e_k: \sum_{k=1}^j \beta_k = \alpha_i, \text{ for } i\ge j+1,\ \beta_i \ge 0\right\} \right)\\
		&= \conv\left( \{\alpha_i e_i\}_{i=1}^j \cup \bigcup_{i=j+1}^n \alpha_i \Delta_j \right).
	\end{align*}
	However, since $\{\alpha_i e_i\}_{i=1}^j$ are the largest,  $\alpha_i \Delta_j  \subseteq \conv(\{\alpha_i e_i\}_{i=1}^j)$, and we conclude that the Newton polytope of $g_\ell$ is the simplex $\tilde{K} = \conv(\{\alpha_i e_i\}_{i=1}^j).$
	
	We have found that the number of solutions is equal to the mixed volume $V(\tilde{K}[j], \Delta_n [n-j]),$ i.e. that
	\[V_n(\tilde{K}[j], \Delta_n [n-j]) = V_n({K}[j], \Delta_n [n-j]). \]
	Since $\tilde{K}$ is a subset of a $j$-dimensional subspace $E = \sp \{e_i\}_{i=1}^j$, we may use a classical calculation regarding mixed volume of lower dimensional sets (attributed to Fedotov, see \cite{burago2013geometric}) to conclude that
	\begin{align*}
	&V_n(\tilde{K}[j], \Delta_n [n-j]) \\
	&= V_n({K}[j], \Delta_n [n-j])\\
	&= \binom{n}{j}^{-1} V_j(\tilde{K}[j])V_{n-j}(P_{E^\perp} \Delta_n [n-j]))\\
	& = \frac{j!(n-j)!}{n!}\vol_j(\tilde{K}) \vol_{n-j}(\Delta_{n-j}) = \frac{\prod_{i=1}^j \alpha_i}{n!}.
	\end{align*}

 \end{proof}

\begin{cor}\label{cor:not-regular-simplex}
		Let $(\alpha_i)_{i=1}^n, (\beta_i)_{i=1}^n$ be two sequences of non-negative numbers, and let $K,T \subset \RR^n$ be given by $K=\conv(0,\alpha_1 e_1,\dots, \alpha_n e_n)$ and $T=\conv(0,\beta_1 e_1,\dots, \beta_n e_n)$. Then, for any $0 \le j \le n$,
		\begin{equation}\label{eq:corollary}
		V_n(K[j],T[n-j])=\frac{1}{n!} \max \left \{\prod_{i\in I} \alpha_i \prod_{j\in I^c} \beta_j : I\subset [n],\ |I|=j\right\}.
		\end{equation}
\end{cor}

\begin{proof}
 	Denote by $A$ the matrix with $(\beta_i^{-1})_{i=1}^n$ on the diagonal, then $AT = \Delta_n$ and $AK =\conv(0,\frac{\alpha_1}{\beta_1} e_1,\dots, \frac{\alpha_n}{\beta_n} e_n)$ , and by Lemma \ref{lem:mixed-of-simpex},
 	\[V(AK[j],AT[n-j])=\frac{\prod_{i\in I} \frac{\alpha_i}{\beta_i}}{n!}\]
 	where $(\frac{\alpha_i}{\beta_i})_{i\in I}$ are the $j$ largest out of $(\frac{\alpha_i}{\beta_i})_{i=1}^n$. We claim that the set $I$ attains the maximum in \eqref{eq:corollary}. Indeed, if not, there is some $i_0 \in I$ and $i_1 \notin I$ such that
 	\[\prod_{i\in I} \alpha_i \prod_{I^c} \beta_j < \alpha_{i_1} \beta_{i_0} \prod_{i\in (I\backslash i_0)} \alpha_i \prod_{j\in (I^c\backslash i_1)} \beta_j.\]
 	Reorganizing this inequality we get $\frac{\alpha_{i_0}}{\beta_{i_0}} < \frac{\alpha_{i_1}}{\beta_{i_1}}$, and so $(\frac{\alpha_i}{\beta_i})_{i\in I}$ are not the $j$ largest, as we may replace $i_0$ with $i_1$, in contradiction to the choice of $I$. 
 	
 	Using the equivariance of mixed volume,
 	\[V(K[j],T[n-j])=\det(A)^{-1} V(AK[j],AT[n-j]) = \frac{\prod_{i\in I} \alpha_i \prod_{j\in I^c}\beta_j}{n!},\]
 	and since $I$ attains the maximum in \eqref{eq:corollary}, the proof is completed.
 	
 \end{proof}
\begin{rems}
	We note that Lemma \ref{lem:mixed-of-simpex} may be used to obtain a calculation of the mixed volumes of a coordinate simplex and a lower dimensional simplex as well.
	Using the $GL_n$ equivariance of mixed volume, the Lemma may be applied to any two ``aligned" simplices, not necessarily with the coordinate directions.
\end{rems}

\bibliographystyle{amsplain}
\addcontentsline{toc}{section}{\refname}\bibliography{AB-inequalities}
\end{document}